\documentclass[11pt]{amsart}
\hoffset= - 0.75 in

\usepackage{amsfonts, amsmath, amssymb, amsthm}
\usepackage{latexsym, graphics}
\newtheorem{thm}{Theorem}

\newtheorem{defn}{Definition}

\begin{document}

\title{Maximal triangulations of a regular prism}
\author{Mike Develin}
\address{Mike Develin, American Institute of Mathematics, 360 Portage 
Ave., Palo Alto, CA 94306-2244, USA}
\date{\today}
\email{develin@post.harvard.edu}

\begin{abstract}

In this paper, we answer two conjectures of De~Loera, Santos, and Takeuchi by demonstrating that the
maximal size of a regular triangulation of a prism over a regular $n$-gon is
$\lceil\frac{n^2+6n-16}{4}\rceil$, and that the maximal size of a regular triangulation of the regular 
$n$-antiprism is $\lfloor\frac{n^2+8n-16}{4}\rfloor$.

\end{abstract}

\maketitle

\section{Introduction}
In a 2001 paper, De~Loera, Santos, and Takeuchi~\cite{DST} showed that the size of the maximal 
triangulation of the prism over an $n$-gon depends on its coordinatization. For a specific 
coordinatization, they 
proved that the maximal triangulation has size between $\lceil\frac{n^2+6n-16}{4}\rceil$ and 
$\frac{n^2+n-6}{2}$. They exhibited a triangulation of a particular coordinatization achieving the 
upper bound, and conjectured that the regular $n$-prism achieved the lower bound.

In this paper, we prove this conjecture. We also prove a related conjecture, namely that the maximal 
triangulation of the regular $n$-antiprism (the convex hull of a regular $n$-gon and its 
image translated one unit in a perpendicular direction and rotated by $\frac{\pi}{n}$) has size 
$\lfloor\frac{n^2+8n-16}{4}\rfloor$. 

\section{Main Proof}

\begin{thm}
Let $P$ be the prism over a regular $n$-gon. Then the maximal triangulation of $P$ has precisely 
$\lceil \frac{n^2+6n-16}{4}\rceil$ 
tetrahedra.
\end{thm}

\begin{proof}
De~Loera, Santos, and Takeuchi~\cite{DST} demonstrated the existence of such a triangulation, given by 
splitting the prism vertically and constructing placing triangulations on both halves.

Our goal is to show that no triangulation can have more than $\lceil \frac{n^2+6n-16}{4}\rceil$
tetrahedra. Certainly there are always $n-2$ tetrahedra joining a triangle from the top $n$-gon to a
vertex on the bottom $n$-gon, and another $n-2$ tetrahedra joining a triangle from the bottom $n$-gon
to a vertex on the top $n$-gon. Therefore, to maximize the number of total tetrahedra, we must
maximize the number of tetrahedra joining an edge from the top $n$-gon to an edge from the bottom
$n$-gon. Given an edge in the top $n$-gon, we define its {\bf link} to be the set of edges in the
bottom $n$-gon to which it is joined to form tetrahedra. We wish to find an upper bound for the total
size of these sets.

\begin{figure} 
\begin{center}\includegraphics{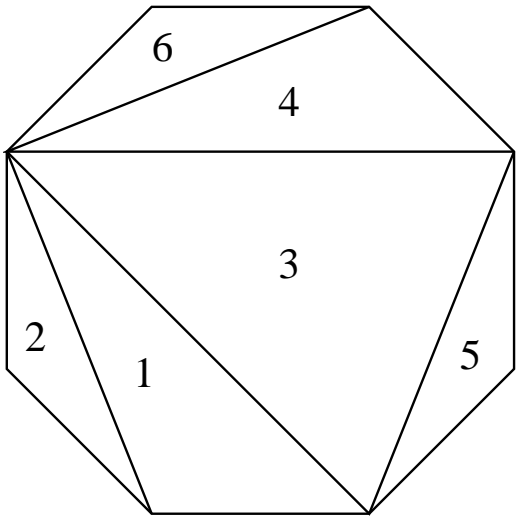}\end{center} 
\caption{\label{lexicog} 
An admissible ordering of the triangles in a triangulation of an 8-gon.}
\end{figure}

Fix a boundary edge of the top $n$-gon; then the dual of the triangulation of the top $n$-gon is a rooted tree on $n-2$ vertices, 
with root equal to the triangle adjacent to that edge~\cite{stancomb}. Let $T_1,\ldots,T_{n-2}$ be an ordering of the triangles in 
this triangulation compatible with the rooted tree, i.e. so that $T_1$ is the root, and if $T_i$ lies above $T_j$ in the tree then 
$i>j$. The resulting numbering is exemplified by Figure~\ref{lexicog}.

Let $x_i$ be the unique
vertex of the bottom $n$-gon is joined to $T_i$. For each edge $e_j$ of the top $n$-gon, define a linear functional $f_j$ as
follows: if $e_j$ is an exterior edge, then let $f_j$ be 0 on $e_j$ and positive on the interior of the $n$-gon, and
if $e_j$ is an interior edge, let $f_j$ be 0 on $e_j$ and positive on the triangle of greater index incident to $e_j$. 

If $e_j$ is an exterior edge incident on some triangle $T_i$, then the link of $e_j$ is a path from $x_i$ to one of the vertices
directly below $e_j$, and it must be monotonely decreasing with respect to $f$ (so that adjacent simplices satisfy the condition
of being on opposite sides of their common facet.) Similarly, if $e_j$ is an interior edge between triangles $T_i$ and $T_k$,
$i<k$, then the link of $e_j$ must be a path from $x_k$ to $x_i$ which is monotonely decreasing with respect to $f_j$.

This motivates the following definition.

\begin{defn}
Let $x$ be a vertex of the (bottom) $n$-gon; let $e_j$ be an edge and $f_j$ the corresponding linear functional. Then 
the {\bf height} $\text{ht}_{e_j}(x)$ of $x$ with respect to $e_j$ is the number of distinct values $f_j(v)$ which 
are (strictly) less than $f_j(x)$, as $v$ ranges over all vertices of the $n$-gon.
\end{defn}

Note that because we are dealing with a regular $n$-gon, most values will be achieved exactly twice, with possible 
singly achieved minimum and maximum. It is precisely this parallelism which allows us to prove this small upper bound 
for the regular $n$-gon.

For an exterior edge $e_j$ adjacent to a triangle $T_i$, the link of $e_j$ is bounded in cardinality by the length of 
the longest monotonely decreasing (with respect to $f_j$) path from $x_i$ to either vertex on $e_j$. This is, by 
definition, at most $\text{ht}_{e_j}(x_i)$. Similarly, for an interior edge $e_j$ adjacent to triangles $T_i$ and 
$T_k$, $i<k$, the size of the link of $e_j$ is at most ${\rm ht}_{e_j}(x_k)-{\rm ht}_{e_j}(x_i)$ (a positive quantity, 
since the link is a monotonely decreasing path with respect to $f_j$.) Therefore, the 
total size of the links of the edges is at most:

\begin{equation}\label{maxlink}
\alpha = 
\sum_{\substack{e_j\,{\rm exterior} \\ e_j \cap T_i\neq \emptyset}} {\rm ht}_{e_j}(x_i) + 
\sum_{\substack{e_j\,{\rm interior} \\ e_j = T_i\cap T_k, i<k}} {\rm ht}_{e_j}(x_k) - {\rm ht}_{e_j}(x_i).
\end{equation}

The ordering we have chosen has the property that for every $i>1$, there exists a unique $\phi(i)<i$ for which $T_{\phi(i)}$
borders $T_i$. Let $e_i$ be the common edge of these two triangles; then by construction, all triangles on the $T_i$ side of $e_i$
have index at least $i$. The region on the $T_i$ side of the diagonal $e_i$ is, of course, a polygon all of whose non-$e_i$ edges
are original edges of the $n$-gon. Let the boundary edges of this polygon (including $e_i$ itself) be denoted by $X_i$, and let
$X_1$ be the set of boundary edges of the original polygon.

We claim that we can 
re-express 
(\ref{maxlink}) as:

\begin{equation}\label{revmaxlink}
\alpha = 
\sum_{e_j\in X_1} {\rm ht}_{e_j}(x_1) + 
\sum_{2\le i\le n-2} (\sum_{e_j\in X_i} {\rm ht}_{e_j}(x_i)-{\rm 
ht}_{e_j}(x_{\phi(i)})).
\end{equation}

Indeed, for every interior edge $e_j$, the function ${\rm ht}_{e_j}$ appears in only one term of (\ref{revmaxlink}), 
namely the $j$-th term; this takes care of the second term of (\ref{maxlink}). For exterior edges, it is a bit more 
complicated; if an exterior edge $e_j$ is incident on some $T_k$, then the $X_i$'s mentioning ${\rm ht}_{e_j}$ have 
their corresponding $T_i$'s forming a path from $T_1$ to $T_k$, and the sum will telescope to yield (when combined 
with the initial contribution ${\rm ht}_{e_j}(x_1)$) precisely ${\rm ht}_{e_j}(x_k)$ as desired. 

\begin{figure}
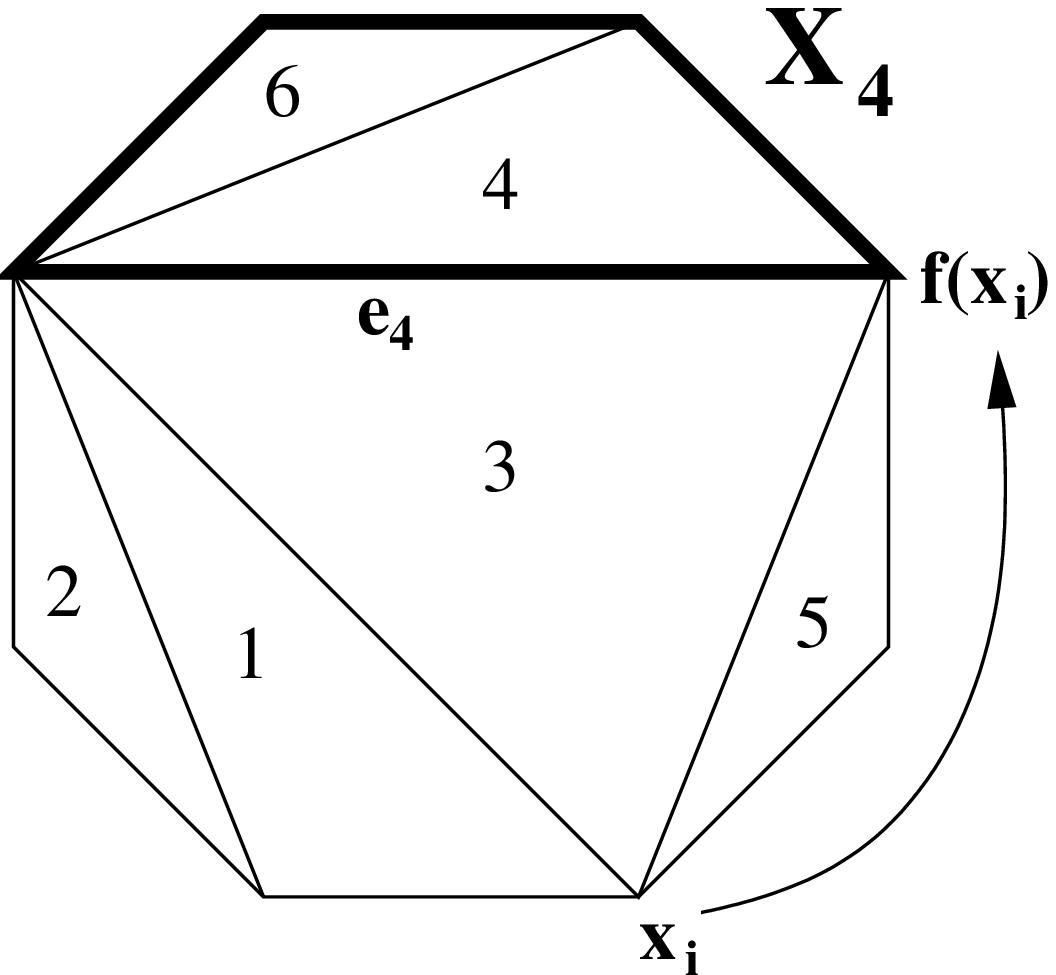
\caption{\label{movinup} 
The vertex $x_i$ must be above $x_{\phi(i)})$ with respect to $e_i$.}
\end{figure}

Now, we evaluate (\ref{revmaxlink}). In particular, we claim that all but the first term must be non-positive. The 
$i$-th term, $i>1$, is equal to $\sum_{e_j\in X_i} {\rm ht}_{e_j}(x_i)-{\rm ht}_{e_j}(x_{\phi(i)})$, or $\sum_{e_j\in 
X_i} {\rm ht}_{e_j}(x_i) - \sum_{e_j\in X_i} {\rm ht}_{e_j}(x_{\phi(i)}))$. We know that ${\rm ht}_{e_i}(x_i) \ge {\rm 
ht}_{e_i}(x_{\phi(i)}))$; see Figure~\ref{movinup} for the current state of affairs, where we rotate the picture if 
necessary to assume that $e_i$ is horizontal for the terminology to follow. Define $g(v) = \sum_{e_j\in X_i} {\rm 
ht}_{e_j}(v)$; it suffices to show that if ${\rm 
ht}_{e_i}(v_1)\ge {\rm ht}_{e_i}(v_2)$, then $g(v_2)\le g(v_1)$. 

\begin{figure}
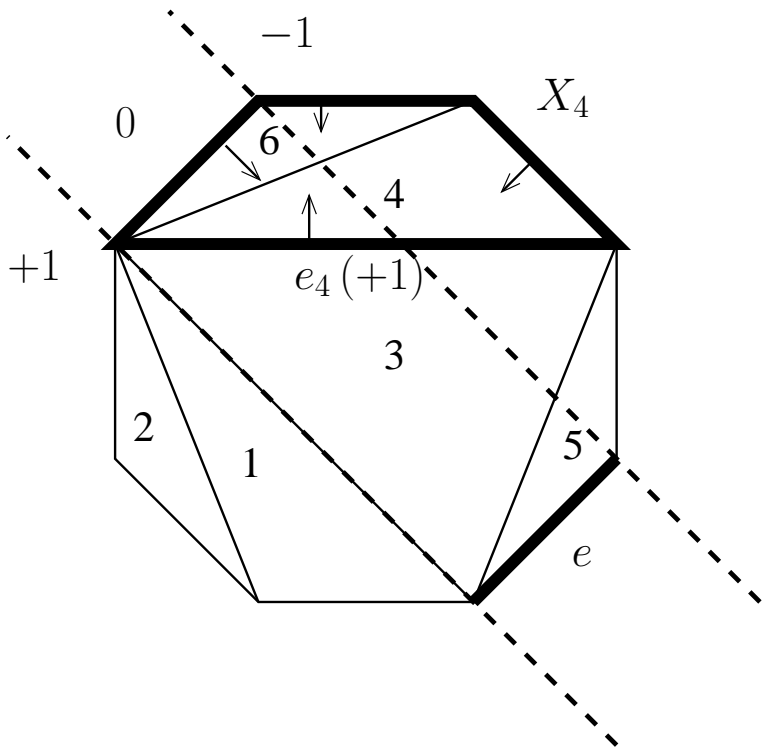
\caption{\label{updown} 
There are fewer edges in $X_4$ (not counting $e_4$) in the +1 region than the -1 region, a straightforward 
consequence of the fact that $v$ is moving up along $e$ relative to $e_4$.}
\end{figure}

The picture is obviously left-right symmetric, so we can assume that $v_1$ and $v_2$ are both on the right side of 
the picture. The vertex $v_2$ lies above $v_1$ by some sequence of steps along the boundary of the $n$-gon; it suffices to show 
that each step up along an edge $e$ does not increase $g(v)$. Indeed, looking at the height difference in each ${\rm 
ht}_{e_j}(v)$, we see that this goes up  by 1 if $e_j$ is on the left side of $e$ around the boundary or equal to 
$e_i$, it goes down by 1 if $e_j$ is on the right side of $e$ around the boundary, and it remains fixed if $e_j$ is 
parallel to $e$ (either equal to it, or its opposite edge in the case where $n$ is even.) Consulting 
Figure~\ref{updown}, which shows the changes in height for boundary edges of the $n$-gon upon moving $v$ 
up along the edge $e$, it is immediate that the number of edges in $X_i$ where the height increases is less than or 
equal to the number of edges where the height decreases, proving the claim.

Consequently, (\ref{revmaxlink}) implies that $\alpha\le\sum_{e_j\in X_1} {\rm ht}_{e_j}(x_1)$. However, this is easy 
to evaluate; it does not matter which vertex $x_1$ is. If $n$ is even, these heights are (going cyclically around the 
$n$-gon) $0,1,2,\ldots,\frac{n}{2}-1,\frac{n}{2}-1,\ldots,2,1,0$, and $\alpha \le (\frac{n}{2}-1)(\frac{n}{2}) = 
\frac{n^2-2n}{4}$; if $n$ is odd, the heights are $0,1,2,\ldots,\frac{n-3}{2}, 
\frac{n-1}{2},\frac{n-3}{2},\ldots,2,1,0$, and we obtain $\alpha \le \frac{n^2-2n+1}{4}$. Recalling that this was a 
bound on the number of simplices formed by joining an edge of the top $n$-gon to an edge of the bottom $n$-gon, and 
adding back the $2n-4$ other simplices yields that the total number of simplices is at most $\frac{n^2+6n-16}{4}$ for 
$n$ even and $\frac{n^2+6n-15}{4}$ for $n$ odd, i.e. the $\lceil \frac{n^2+6n+16}{4}\rceil$ bound in question.
\end{proof}

The same proof technique works to prove the conjectured bound of De Loera, Santos, and Takeuchi on maximal 
triangulations of the regular $n$-antiprism, defined to be the convex hull of a regular $n$-gon and its image one 
unit away in the perpendicular direction and rotated by $\frac{\pi}{n}$.

\begin{thm}
The maximum number of simplices in a triangulation of the regular $n$-antiprism is at most $\lfloor 
\frac{n^2+8n-16}{4}\rfloor$.
\end{thm}

\begin{proof}
De Loera, Santos, and Takeuchi~\cite{DST} demonstrated a triangulation of the regular $n$-antiprism containing this 
number of simplices. By an identical argument to the above one, we can bound this number above by $\sum_{e_j\in X_1} 
{\rm ht}_{e_j}(x)$, where $x$ is any vertex in the bottom $n$-gon and $X_1$ is the set of boundary edges of the 
upper $n$-gon. Computing these heights
yields (in cyclic order) $0,1,\ldots,\frac{n-1}{2},\frac{n-1}{2},\ldots,1$ for $n$ odd, and 
$0,1,\ldots,\frac{n}{2}-1, \frac{n}{2},\frac{n}{2}-1,\ldots,1$ for $n$ even. This sum is $\frac{n^2-1}{4}$ for $n$ 
odd 
and $\frac{n^2}{4}$ for $n$ odd. Adding in the $2n-4$ simplices formed by joining a triangle in the upper $n$-gon to 
a vertex in the lower $n$-gon or vice versa yields the desired bound.
\end{proof}

\section*{Acknowledgements}
The author was supported by an American Institute of Mathematics Postdoctoral Fellowship. I would also like to thank Francisco 
Santos for helpful comments on an early draft.


\begin{thebibliography}{99} 
\bibitem{DST}
J. De Loera, F. Santos, F. Takeuchi, ``Extremal properties for dissections of convex 3-polytopes'', \textit{Siam J. 
Discrete Math.} \textbf{14} (2001), 143--161.

\bibitem{stancomb}
R. Stanley, \textit{Enumerative Combinatorics, Volume II}, Cambridge University Press, 1999. 

\end{thebibliography}
\end{document}